\documentclass{article}[13pt]
\usepackage{amssymb, latexsym, amsmath, amsthm, hyperref}
\newtheorem{theorem}{Theorem}
\newtheorem{lemma}[theorem]{Lemma}

\newtheorem{corollary}[theorem]{Corollary}
\newtheorem{proposition}[theorem]{Proposition}

\begin{document}

\centerline{\textbf{\Large{Observations on some classes of operators on C(K,X)}}}

\vspace{.1in}
\centerline{Ioana Ghenciu}

\vspace{.1in}
\noindent
%Mathematics Department, 
University of Wisconsin-River Falls, 
 River Falls, Wisconsin, 54022, USA
	
	\noindent
	email: ioana.ghenciu@uwrf.edu 
	
	\centerline{Roxana  Popescu}
	\vspace{.1in}
	\noindent
	University of Pittsburgh, Pittsburgh, PA, 15260, USA.

	\noindent
	email: rop42@pitt.edu

\begin{abstract}

Suppose  $X$ and $Y$ are Banach spaces,  $K$ is a compact Hausdorff space,  $\Sigma$ is the $\sigma$-algebra of Borel subsets of $K$,  $C(K,X)$ is the Banach space of all continuous $X$-valued functions (with the supremum norm), and $T:C(K,X)\to Y$ is a strongly bounded operator with representing measure $m:\Sigma \to L(X,Y)$. 

%$T$ is a strongly bounded operator  and
We  show that if  $\hat{T}: B(K, X) \to Y$ is its extension,  then $T$ is weak Dunford-Pettis (resp. weak$^*$ Dunford-Pettis, weak $p$-convergent, weak$^*$ $p$-convergent)  if and only if  $\hat{T}$ has the same property. 
%is weak Dunford-Pettis (resp. weak$^*$ Dunford-Pettis, weak $p$-convergent, weak$^*$ $p$-convergent,  $1\le p< \infty$). 

%We give characterizations of strongly bounded  limited completely continuous (resp. limited $p$-convergent) operators in terms of their representing measures. 

We prove that if $T:C(K,X)\to Y$   is strongly bounded limited completely continuous (resp.  limited $p$-convergent), then  $m(A):X\to Y$ is limited completely continuous (resp.  limited $p$-convergent) for each $A\in \Sigma$. We also prove that  the above implications become equivalences when  $K$ is a dispersed compact Hausdorff space. 

% if $K$ is a dispersed compact Hausdorff space  and $T$ is a strongly bounded operator, then $T$ is  limited completely continuous (resp.  limited $p$-convergent) whenever $m(A):X\to Y$ is limited completely continuous (resp.  limited $p$-convergent), for each $A\in \Sigma$.

%and  $m\leftrightarrow T: C(K,X)\to Y$ is a strongly bounded operator such that  for each $A\in \Sigma$, $m(A):X\to Y$ is limited completely continuous (resp.  limited $p$-convergent), then $T$  is limited completely continuous (resp.  limited $p$-convergent).

\end{abstract}

\noindent Key words and phrases: weak Dunford-Pettis operators, weak$^*$ Dunford-Pettis operators, weak $p$-convergent operators,  weak$^*$ $p$-convergent operators, limited completely continuous operators,  limited $p$-convergent operators\\
%pseudo weakly compact operators, $p$-Dunford-Pettis completely continuous operators, 
[.1in]
Subject classification:  46B20, 46B25, 46B28\\
%[.1in]
Mathematical discipline: Functional Analysis

\section{Introduction}

Suppose $K$ is a compact Hausdorff space, $X$ and $Y$ are Banach spaces, $C(K,X)$ is the Banach space of all continuous $X$-valued functions (with the supremum norm), and  $\Sigma$ is the $\sigma$-algebra of Borel subsets of $K$.   
The topological dual of $C(K, X)$ can be identified with the space $rcabv(\Sigma, X^*)=M(\Sigma, X^*)$  of all $X^*$-valued, regular countably additive Borel measures on $K$ of bounded variation, endowed with the variation norm.
%(or $M(K, X^*)$)

Every continuous linear operator $T: C(K,X) \to Y $ may be represented by a vector measure $m: \Sigma \to L(X,Y^{**})$ of finite semi-variation (\cite{BrL},  \cite[p.182]{DU}) %\cite{D},
such that
$$ T(f)= \int_K f \, dm \, \, \, , \, f\in C(K,X), \, \,  \,\, \|T\|=\tilde{m}(K),  $$
%\text{and}
and $T^*(y^*)=m_{y^*}$, $y^*\in Y^*$, where $\tilde{m}$ denotes the semivariation of $m$. 
For each  $y^*\in Y^*$, the vector measure $m_{y^*}=y^*m:\Sigma\to X^*$ defined by  $\langle m_{y^*}(A), x\rangle =\langle  m(A)(x), y^* \rangle $, $A\in \Sigma$, $x\in X$, is a regular countably additive measure of bounded variation. We denote the correspondence $m \leftrightarrow T$.
If we denote by $|y^*m |$ the variation of the measure $y^*m$, then for $E\in \Sigma$, the semi-variation $\tilde{m}(E)$ is given by 
$$\tilde{m}(E)=\sup\{ |y^*m|(E): y^*\in Y^*, \|y^*\|\le 1\}. $$

%	%(\cite{BrL}, lemma 3.1). 

We note that for $f\in C(K,X)$, $\int_K f \, dm \in Y$ even if $m$ is not $L(X,Y)$-valued.
A representing measure $m$ is called \emph{strongly bounded}  if $(\tilde{m}(A_i))\to 0$ for every decreasing  sequence $(A_i)\to \emptyset $ in $\Sigma$,
and  an operator $m \leftrightarrow T : C(K, X)\to Y$ is called strongly bounded if $m$ is strongly bounded  \cite{BrL}.
   By Theorem 4.4 of \cite{BrL}, a strongly bounded representing measure takes its values in $L(X,Y)$. 
	%We note that $m$ is strongly bounded if and only if $\{|y^*m|: y^*\in Y^*, \|y^*\|\le 1\}$ is uniformly countably additive (\cite[Lemma 3.1]{BrL}).  Consequently, if 
	If $m$ is a strongly bounded representing measure, then there is a non-negative regular Borel measure $\lambda$ so that $\tilde{m}(A)\to 0$ as $\lambda(A) \to 0$ (\cite[Lemma 3.1]{BrL} and the proof of \cite[Theorem 4, p. 11]{DU}). We call the measure $\lambda$ a \emph{control measure} for $m$.  
	%(see \cite[Lemma 3.1]{BrL}). from Strictly Singular and strictly cosingular operators etc, Bator and Lewis

Let $\chi_A$ denote the characteristic function of a set $A$, and $B(\Sigma, X)$  denote the space of all bounded, $\Sigma$-measurable functions on $K$ with separable range   in $X$ and the sup norm. Certainly $C(K,X)$ is contained isometrically in $B(\Sigma, X)$. Further, $B(\Sigma, X)$ embeds isometrically in $C(K,X)^{**}$; e.g., see \cite{BrL}. The reader should note that if $m \leftrightarrow T$, then $m(A)x=T^{**}(\chi_A x)$, for each $A\in \Sigma$, $x\in X$. 
If $f\in B(\Sigma,X)$, then $f$ is the uniform limit of $X$-valued simple functions, $\int_K f \,dm$ is well-defined and defines an extension $\hat{T}$ of $T$; e.g., see \cite{D}.
Theorem 2 of \cite{BB} shows that $\hat{T}$ maps $B(\Sigma, X)$ into $Y$ if and only if the representing measure
$m$ of $T$ is $L(X,Y)$-valued. If $T:C(K,X) \to Y$ is strongly bounded, then $m$ is $L(X,Y)$-valued  \cite{BrL}, and thus $\hat{T}: B(\Sigma, X)\to Y$.
Since  $\hat{T}$ is the restriction to $B(\Sigma, X)$ of the operator $T^{**}$, it is clear that an operator $T: C(K,X)\to Y$ is compact (resp. weakly compact) if and only if its extension $\hat{T}: B(\Sigma,X) \to Y$ is compact (resp. weakly compact).

%\cite{GhActa},
%weakly compact, compact, 
Several authors have found the study of $\hat{T}$ to be quite helpful. We mention the work of \cite{BB},  \cite{BCCam}, \cite{BPNach},  \cite{GLG},  \cite{CS},  \cite{GLBCl},  \cite{MA}, \cite{GhActa}, and \cite{GPQM}. In these papers it has been proved that if $m$ is strongly bounded, then $T: C(K,X)\to Y$ is 
Dunford-Pettis, Dieudonn\'{e}, unconditionally converging,  strictly singular, strictly cosingular, weakly precompact,  weakly $p$-compact,  limited, pseudo weakly compact, Dunford-Pettis $p$-convergent if and only if its extension $\hat{T}: B(\Sigma,X)\to Y $ has the same property.  
Further, $T$ has a weakly precompact, Dunford-Pettis, unconditionally converging, 
$p$-convergent,
Dunford-Pettis $p$-convergent, pseudo weakly compact, limited completely continuous, limited $p$-convergent adjoint if and only if its extension $\hat{T}: B(\Sigma,X)\to Y $ has the same property.  
%$p$-convergent,

%$1<p<\infty$
We show that if $T: C(K,X)\to Y$ is a  strongly bounded  operator, $\hat{T}: B(\Sigma, X) \to Y$ is its extension, 
then $T$ is weak Dunford-Pettis (resp. weak$^*$ Dunford-Pettis, weak $p$-convergent, weak$^*$ $p$-convergent, $1\le p<\infty$)  if and only if  $\hat{T}$ is weak Dunford-Pettis (resp. weak$^*$ Dunford-Pettis, weak $p$-convergent, weak$^*$ $p$-convergent, $1\le p<\infty$). 
It is shown that if  $m\leftrightarrow T: C(K,X)\to Y$ is strongly bounded and  $T$ is weak Dunford-Pettis (resp. weak$^*$ Dunford-Pettis, weak $p$-convergent, weak$^*$ $p$-convergent), then $m(A):X\to Y$ is weak Dunford-Pettis (resp. weak$^*$ Dunford-Pettis, weak $p$-convergent, weak$^*$ $p$-convergent) for each $A\in \Sigma$. If we additionally assume that  $K$ is a dispersed compact Hausdorff space, then
the above implications become equivalences.  
% a strongly bounded operator  $m\leftrightarrow T: C(K,X)\to Y$  is  weak Dunford-Pettis (resp. weak$^*$ Dunford-Pettis, weak $p$-convergent, weak$^*$ $p$-convergent) if and only if for each $A\in \Sigma$, $m(A):X\to Y$ is   weak Dunford-Pettis (resp. weak$^*$ Dunford-Pettis, weak $p$-convergent, weak$^*$ $p$-convergent). We use the same ideas and techniques as in \cite{BRS}, \cite{BCCam}, and \cite{CS}. 

We give characterizations of strongly bounded pseudo weakly compact (resp. Dunford-Pettis $p$-convergent,  limited completely continuous limited $p$-convergent) operators in terms of their representing measures. 
% Th 15 and Remark, Corollary 18,  and Th 20

\section{Definitions and Notation} 
\vspace{0.1in}

%, $E$ and $F$ 
Throughout this paper, $X$ and $Y$ will denote Banach spaces. The unit ball of $X$ will be denoted by $B_X$,   
and  $X^*$ will denote the continuous linear dual of $X$. 
 %$X^*$ will denote the continuous linear dual of $X$, and the closed linear span of a sequence $(x_n)$ in $X$ will be denoted by $[x_n]$. 
%The space $X$ embeds in $Y$ (in symbols $X\hookrightarrow Y$) if $X$ is isomorphic to a closed subspace of $Y$. 
An operator $T: X \to Y$ will be a continuous and linear function. 
%We denote by $L(X,Y)$ (resp. $W(X,Y)$) the space of all bounded linear operators (resp. weakly compact)  from $X$ to $Y$ 
%The set of all operators, resp. weakly compact operators from $X$ to $Y$  will be denoted by $L(X,Y)$, resp.   $W(X,Y)$. 
The set of all operators from $X$ to $Y$  will be denoted by $L(X,Y)$.

An operator $T:X\to Y$ is called \emph{completely continuous (or Dunford-Pettis)} if $T$ maps weakly convergent sequences to norm convergent sequences. 

A Banach space $X$ has the \emph{Dunford-Pettis property (DPP)} if every weakly compact operator $T:X \to Y$ is completely continuous, for any Banach space $Y$. 
%(or Dunford-Pettis)

A bounded subset $A$ of $X$ is called a \emph{limited (resp. Dunford-Pettis or $DP$}) subset  of $X$ if each $w^*$-null (resp. weakly null) sequence $(x_n^*)$ in $X^*$ tends to $0$ uniformly on $A$; i.e.
$$\sup_{x \in A} |x_n^* (x) | \to  0.$$

A subset $S$ of $X$ is said to be \emph{weakly precompact} (or \emph{conditionally weakly compact}) provided that every  sequence from $S$ has a weakly Cauchy subsequence.
%	Every DP set is weakly precompact, e.g., see \cite[p. 377]{HR},  \cite{KA}. 
An operator   $T:X\to Y$ is called \emph{ weakly precompact } if $T(B_X)$ is weakly precompact.

A series $\sum x_n$ in $X$ is said to be \emph{weakly unconditionally convergent (wuc)} if for every $x^*\in X^*$, the series $\sum |x^*(x_n)|$ is convergent. An operator $T:X\to Y$ is  \emph{unconditionally converging} if it maps weakly unconditionally convergent series to unconditionally convergent ones.

% $p^*$
For $1\le p<\infty$, $p'$ denotes the conjugate  of $p$. If $p=1$, then  $c_0$ plays the role of $\ell_{p'}$.
The unit vector basis of $\ell_p$ will be denoted by $(e_n)$. 

%\le 
Let $1\le p<\infty$. A sequence $(x_n)$ in $X$  is called  \emph{weakly} $p$-\emph{summable}  if $(x^*(x_n))\in \ell_p$ for each $x^*\in X^*$ \cite[p. 32]{DJT}, \cite[p. 134]{Ryan}. 
Let $\ell_p^w(X)$ denote the  set of all weakly $p$-summable sequences in $X$. The space $\ell_p^w(X)$ is a Banach space with the norm 
$$ \|(x_n)\|_p^w=\sup \{( \sum_{n=1}^\infty |\langle x^*, x_n \rangle |^p )^{1/p}: x^*\in B_{X^*}\}$$
%\begin{align*} \|(x_n)\|_p^w& =\sup \{( \sum_{n=1}^\infty |\langle x^*, x_n \rangle |^p )^{1/p}: x^*\in B_{X^*}\}\\
%&=\sup \{\|\sum_{n=1}^\infty a_n x_n \|: a=(a_n)\in B_{\ell_{p'}} \}, \qquad \text{ when} \quad 1<p<\infty,   
%\end{align*}
%and 
%\begin{align*} \|(x_n)\|_1^w& =\sup \{\sum_{n=1}^\infty |\langle x^*, x_n \rangle | : x^*\in B_{X^*}\}\\
%&=\sup \{\|\sum_{n=1}^\infty a_n x_n \|: a=(a_n)\in B_{c_0} \}, \qquad \text{ when} \quad p=1.    
%\end{align*}

%If $p=\infty$, then $\ell_\infty(X)=\ell_\infty^w(X)$ \cite[p. 33]{DJT}; if $(x_n)$ is a bounded sequence in $X$, then 
%$$\|(x_n) \|_\infty^w= \sup_n \|x_n\|=\|(x_n)\|_\infty. $$

Let $c_0^w(X)$ be the space of weakly null sequences in $X$.
% This is a  Banach space with the norm 
%$$\|(x_n)\|_{c_0^w}=\sup_{\|x^*\|\le 1} \|(\langle x^*, x_n \rangle)\|_{c_0}, $$
%and $c_0^w(X)\simeq W(\ell_1, X)$ \cite{FZ2}. 
%, \cite[Proposition 1.2.5]{FZ3}. 
For $p=\infty$, we consider the space $c_0^w(X)$ instead of $\ell_{\infty}^w(X)$. 
%=\ell_{\infty}(X)

If $p<q$, then  $\ell_p^w(X)\subseteq \ell_q^w(X)$. Further, the unit vector basis of $\ell_{p'}$ is  weakly $p$-summable for all $1<p<\infty$. The weakly  $1$-summable sequences are precisely  the weakly unconditionally convergent  series
and the weakly $\infty$-summable sequences are precisely the weakly null sequences. 

We recall the following  isometries: 
 $L(\ell_{p'},X) \simeq \ell_p^w(X)$ for $1<p<\infty$; $ L(c_0,X)\simeq \ell_p^w(X)  $ for $p=1$; 
$T\to (T(e_n))$   \cite[Proposition 2.2, p. 36]{DJT}, \cite[p. 134]{Ryan}.

Let $1\le p\le \infty$. An operator $T:X\to Y$ is called \emph{p-convergent} if $T$ maps weakly $p$-summable sequences into norm null sequences \cite{CS}. The set of all $p$-convergent operators is denoted by $C_p(X,Y)$. 

The $1$-convergent operators are precisely the unconditionally converging operators and the $\infty$-convergent operators are precisely the completely continuous operators. If $p<q$, then $C_q(X,Y)\subseteq C_p(X,Y)$. 

	%A sequence $(x_n)$ in $X$ is called \emph{weakly p-convergent} to $x\in X$ if the sequence $(x_{n}-x)$ is weakly $p$-summable \cite{CS}. The weakly $\infty$-convergent sequences are precisely the weakly convergent sequences. 

%	Let $1\le p\le \infty$. A bounded subset $K$ of $X$ is \emph{relatively weakly $p$-compact} (resp. \emph{weakly p-compact}) if every sequence in $K$ has a weakly $p$-convergent subsequence with limit in $X$ (resp. in $K$). 
	
	%An operator $T:X\to Y$ is \emph{weakly $p$-compact} if $T(B_X)$ is relatively weakly $p$-compact \cite{CS}. The set of all weakly $p$-compact operators $T:X\to Y$ will be denoted by $W_p(X,Y)$. 

%The weakly $\infty$-compact operators are precisely the weakly compact operators. If $p<q$, then $W_p(X, Y)\subset W_q(X, Y)$ \cite{CS}.

%(resp. $X\in W_p$) (resp. $id(X)\in W_p(X,X)$)

A Banach space $X\in C_p$  if $id(X)\in C_p(X,X)$, where $id(X)$ is the identity map on $X$ \cite{CS}.

A topological space $S$ is called   \emph{dispersed   (or scattered)} if every nonempty closed subset of $S$ has an isolated point. A compact Hausdorff space $K$ is dispersed   if and only if $\ell_1 \not \hookrightarrow C(K)$ \cite[Main theorem]{PS}.

%A sequence $(x_n)$ in $X$ is called \emph{weakly $p$-Cauchy} if $(x_{n_k}-x_{m_k})$ is weakly $p$-summable for any increasing sequences $(n_k)$ and $(m_k)$ of positive integers. 

	%Every weakly $p$-convergent sequence is weakly $p$-Cauchy, and the weakly $\infty$-Cauchy sequences are precisely the weakly Cauchy sequences.  

 %	Let $1\le p\le \infty$. A subset $S$ of $X$ is called \emph{weakly $p$-precompact} if  every  sequence from $S$ has a weakly $p$-Cauchy subsequence.  An operator   $T:X\to Y$ is called \emph{ weakly $p$-precompact  } if $T(B_X)$ is weakly $p$-precompact.
	
	%\cite{IGQM}. 
	%The weakly $\infty$-precompact sets are precisely the weakly precompact sets.
	%A weakly-p-precompact set is also called conditionally weakly p-compact set \cite{CCDL}. 

\section{Operators and Measures}

An operator $T:X\to Y$ is called \emph{weak Dunford-Pettis} \cite[p. 349]{AB} if $\langle T(x_n), y_n^*\rangle \to 0$,  whenever $(x_n)$ is a weakly null sequence in $X$ and $(y_n^*)$ is a  weakly null sequence in $Y^*$.

An operator $T:X\to Y$ is called \emph{weak$^*$ Dunford-Pettis} \cite{KHMB} if $\langle T(x_n), y_n^*\rangle \to 0$,  whenever $(x_n)$ is a weakly null sequence in $X$ and $(y_n^*)$ is a  $w^*$-null sequence in $Y^*$.

Let $1\le p<\infty$. An operator $T:X\to Y$ is called \emph{weak $p$-convergent}  if $\langle y_n^*,T(x_n) \rangle\to 0$ whenever $(x_n)$ is weakly $p$-summable in $X$ and $(y_n^*)$ is weakly null in $Y^*$ \cite{FZ2}.

Let $1\le p<\infty$. An operator $T:X\to Y$ is called \emph{weak$^*$ $p$-convergent}  if $\langle y_n^*,T(x_n) \rangle\to 0$ whenever $(x_n)$ is weakly $p$-summable in $X$ and $(y_n^*)$ is $w^*$-null in $Y^*$ \cite{FZ2}. 

These classes of operators were studied in \cite{AB}, \cite{KHMB}, \cite{FZ2}, and \cite{GRen}. 
In the sequel we will use the following results. 

\begin{lemma} \label{le1}
(\cite[Theorem 2]{BCCam}, \cite[Proposition 1.7.1, p. 36]{CemMen}) A bounded sequence $(f_n)$ in $C(K,X)$ is weakly null if and only if $(f_n(t))$ is weakly null for every $t\in K$. 

%(ii) (\cite[Corollary 1.7.1, p. 36]{CemMen}) A bounded sequence $(f_n)$ in $C(K,X)$ is weakly Cauchy if and only if $(f_n(t))$ is weakly Cauchy for every $t\in K$. 
\end{lemma}

%\begin{proof}(iii) 
%The proof is similar to the proof of \cite[Proposition 1.7.1, p. 36]{CemMen}.
 %Suppose $(f_n)$  is weakly $p$-summable in $C(K, X)$. Let $t\in K$. For each $x^*\in X^*$, let  $\delta_t(\cdot) x^* \in C(K,X)^*$,
 %$$ \langle \delta_t(\cdot) x^*, f\rangle=\int_K f d(\delta_t(\cdot) x^*)= \langle x^*, f(t)\rangle, \qquad f\in C(K, X).$$
%Since $(f_n)$  is weakly $p$-summable, 
 %$$\sum_{n=1}^\infty |\langle \delta_t(\cdot) x^*, f_n\rangle|^p =\sum_{n=1}^\infty | \langle x^*, f_n(t)\rangle|^p<\infty. $$
%Therefore  $(f_n(t))$ is weakly  $p$-summable. 
%\end{proof}

%Lemma \ref{lemma1}

%closed
\begin{lemma} \label{lemma2}
Let $X$ be a Banach space and  $Y$ be a  subspace of $X$. 
%Suppose there is an  operator $T:Y^*\to X^*$ such that $T(y^*)|_Y=y^*$ for every $y^*\in Y^*$. 
Let $(y_n)$ be a sequence in $ Y$. 

(i) If  $(y_n)$ is  weakly null in $X$, then it is  weakly null in $Y$. 

(ii) Let $1\le p<\infty$. If $(y_n)$ is   weakly $p$-summable in $X$, then it is  weakly $p$-summable in $Y$. 
%\end{observation}
\end{lemma}

%%Let $T:Y^*\to X^*$ be an extension operator as in the hypothesis. 

\begin{proof}
We only prove (i); the proof of (ii) is similar. 

(i)  Suppose $(y_n)$ is weakly null  in $X$.
Let $y^*\in Y^*$. Let $x^*\in X^*$ be a Hahn-Banach extension of $y^*$. Then  $ y^*(y_n)= x^*(y_n)\to 0$, since $(y_n)$ is weakly null in $X$. Thus $(y_n)$ is weakly null in $Y$.
\end{proof}

%Suppose that $T: C(K,X)\to Y$ is an operator and   $\hat{T}: B(\Sigma, X) \to Y^{**}$ is its extension to $B(\Sigma,X)$. If $m\leftrightarrow T:C(K,X) \to Y$ is strongly bounded, then $m$ is $L(X,Y)$-valued and $\hat{T}$ maps $B(\Sigma, X)$ into $Y$.

%\vspace{0.1in}
Let  $B_0$ denote the unit ball of $C(K,X)$, and $B$  denote the unit ball of $B(\Sigma, X)$. 

%$1< p< \infty$.
\begin{theorem} \label{T1}
Let $1\le p< \infty$. Suppose that $T: C(K,X)\to Y$ is a  strongly bounded  operator, $\hat{T}: B(\Sigma, X) \to Y$ is its extension. 
Then $T$ is weak Dunford-Pettis (resp. weak$^*$ Dunford-Pettis, weak $p$-convergent, weak$^*$ $p$-convergent)  if and only if  $\hat{T}$ is weak Dunford-Pettis (resp. weak$^*$ Dunford-Pettis, weak $p$-convergent, weak$^*$ $p$-convergent). 
\end{theorem}

\begin{proof} 
We only prove the result for  weak  Dunford-Pettis (resp. weak $p$-convergent) operators. The proof for the other cases is similar.

 Suppose that $T: C(K,X)\to Y$ is weak Dunford-Pettis (resp. weak $p$-convergent) and  $\hat{T}$ is not weak Dunford-Pettis (resp. weak $p$-convergent). 
 Let $(y_n^*)$  be   weakly null in $Y^*$ and $(f_n)$ be a weakly null (resp.  weakly $p$-summable) sequence in  $B(\Sigma, X)$ such that
$ |\langle y_n^*, \hat{T}(f_n) \rangle |  >\epsilon $ for all $n$.   Without loss of generality assume $\|y_n^*\|\le 1 $ for all $n$ and $(f_n)$ is in $B$. 

Using the existence of a control measure for $m$ and  Lusin's theorem, we can find a compact subset $K_0$ of $K$ such that  $\tilde{m} (K\setminus K_0)< \epsilon/4$ and $g_n= f_n|_{K_0}$ is continuous for each $n\in \mathbb{N}$.
 Let $r:B(\Sigma,X)\to B(\Sigma|_{K_0}, X)$ be the restriction operator. Since $(f_n)$ is weakly null (resp. weakly $p$-summable) in $B(\Sigma,X)$, $(g_n)=(r(f_n))$ is weakly null (resp. weakly $p$-summable) in $B(\Sigma|_{K_0},X)$. Every measure $\mu \in C(K_0, X)^*$ can be naturally extended to an element of $B(\Sigma|_{K_0},X)^*$. The sequence $(g_n)$ is weakly null (resp. weakly $p$-summable) in $C(K_0, X)$ by Lemma \ref{le1} (resp. Lemma \ref{lemma2}). 
%Lemma \ref{lemma1}

 Let $H=[g_n]$ be the closed linear subspace spanned by $(g_n)$ in $C(K_0,X)$ and let $S:H\to C(K,X)$  be the isometric extension operator given by   \cite[Theorem 1]{BCCam}.
If $h_n= S(g_n)$ for each $n\in \mathbb{N}$, then $(h_n)$ is  weakly null (resp. weakly $p$-summable) in $C(K, X)$, $(h_n)$ is in the unit ball of $C(K,X)$, and
\begin{align*}
|\langle y_n^*, T(h_n) \rangle | & \ge |\langle y_n^*, \int_{K_0} h_n \,dm  \rangle| - |\langle y_n^*, \int_{K\setminus K_0} h_n \,dm  \rangle|\\
& \ge |\langle y_n^*, \int_{K_0} f_n \,dm  \rangle| - \epsilon/4 \\
& \ge | \langle y_n^*, \int_{K} f_n \,dm  \rangle| - |\langle y_n^*, \int_{K\setminus K_0} f_n \,dm  \rangle|- \epsilon/4 \\
& \ge |\langle y_n^*, \hat{T} (f_n)   \rangle |-  \epsilon/4-\epsilon/4  > \epsilon/2.
\end{align*}
This is a contradiction, since $T$ is a weak Dunford-Pettis (resp. weak $p$-convergent) operator. 
\end{proof}

%$1< p< \infty$.
\begin{corollary} \label{mwDP}
Let $1\le p< \infty$. Suppose that  $m\leftrightarrow T: C(K,X)\to Y$ is a strongly bounded operator. 
If $T$ is weak Dunford-Pettis (resp. weak$^*$ Dunford-Pettis, weak $p$-convergent, weak$^*$ $p$-convergent), then $m(A):X\to Y$ is weak Dunford-Pettis (resp. weak$^*$ Dunford-Pettis, weak $p$-convergent, weak$^*$ $p$-convergent) for each $A\in \Sigma$. 
\end{corollary}

%(resp. weak $p$-convergent) (resp. weak $p$-convergent)

\begin{proof}
If $A\in \Sigma$, $A\neq \emptyset$, define $\theta_A: X\to B(\Sigma,X)$ by $\theta_A(x)= \chi_A x$. Then $\theta_A$ is an isometric embedding of $X$ into $B(\Sigma,X)$ and $ \hat{T}\theta_A= m(A)$. If $T$ is weak Dunford-Pettis, then $\hat{T}$ and thus $m(A)$ is weak Dunford-Pettis, by Theorem \ref{T1}.
\end{proof}

In the next results we will need the following well-known result (\cite{D}, Sec 13, Theorem 5). 

\begin{theorem} \label{Din}
Let $\lambda$ be a positive Radon measure on $K$. If $m\in rcabv(\Sigma, X^*)$ is $\lambda$-continuous, then there exists a function $g:K\to X^*$ such that 

i) $\langle x, g \rangle $ is a $\lambda$-integrable function for every $x\in X$.

ii) For every $x\in X$ and $A\in \Sigma$, 
$$ \langle x, m(A) \rangle =\int_A \langle x, g \rangle \; d\lambda $$

iii) $|g|$ is $\lambda$-integrable and 
$$ \tilde{m}(A)  =\int_A| g |\; d\lambda, $$
for every $A\in \Sigma$, where $|g|(t)=\|g(t)\|$, for $t\in K$. 
\end{theorem}

%m\leftrightarrow
By \cite[Lemma 1.3]{BRS} the operator $ T:C(K, X)\to Y$ with representing measure $m$ is strongly bounded if and only if there exists a positive Radon measure $\lambda$ on $K$ such that

(i) $T^*(y^*)$ is $\lambda$-continuous for every $y^*\in Y^*$, and 

(ii) If $g_{y^*}$ is the element corresponding to $T^*(y^*)$ by Theorem \ref{Din}, then the set $\{|g_{y^*}|:y^*\in B_{Y^*}\}$ is relatively weakly compact in $L^1(\lambda)$. 
 
In this case, $\lambda$ is a control measure for $m$.

%Suppose $m\leftrightarrow T:C(K, X)\to Y$ is strongley bounded,  $\lambda$, $y^*\in Y^*$, and $g_{y^*}$ is as above. Recall that $T^*(y^*)=m_{y^*}\in rcabv(\Sigma, X^*)$, $y^*\in Y^*$. By Theorem \ref{Din}, 
%$\langle T(f), y^* \rangle =\int_K \langle g_{y^*}, f \rangle \;d\lambda$,  $y^*\in Y^*$, for all $f\in C(K,X)$. 
%for every $X$-valued $\Sigma$ simple function $f$ defined  on $K$, and thus for every $f\in C(K,X)$.

\vspace{0.1in}

The following result is motivated by  \cite[Theorem 1.4]{BRS}. 
%similar to

\begin{theorem} \label{BomRS1}
Let $m\leftrightarrow T:C(K, X)\to Y$ be a strongly bounded operator whose representing measure $m$ has a control measure $\lambda$. 
The following are equivalent:

a) $T$ is weak Dunford-Pettis (resp.  weak$^*$ Dunford-Pettis). 
%, weak $p$-convergent, weak$^*$ $p$-convergent

% $Y^*$
b) For every weakly null sequence $(f_n)$ in $C(K, X)$ and every weakly null (resp. $w^*$-null) sequence $(y_n^*)$ in $B_{Y^*}$, we have 
$$\lim_{n\to \infty} \int_K |\langle f_n, g_n\rangle|\; d\lambda=0,  $$
where $g_n$ is the function corresponding to $T^*(y_n^*)$ by Theorem \ref{Din}. 
\end{theorem}

\begin{proof} 
%We give the proof for pseudo weakly compact operators; the other cases are similar. 
$a)\Rightarrow b)$ Suppose $T$ is weak Dunford-Pettis.  Let $(f_n)$ be  a weakly null sequence in $C(K, X)$ and $(y_n^*)$ be a weakly null sequence in $B_{Y^*}$.
% $Y^*$
Let $\phi_n$ be a scalar continuous function on $K$ such that $\|\phi_n\|\le 1 $ and 
$$ \int_K |\langle f_n, g_n\rangle|\; d\lambda \le  \int_K \phi_n \langle f_n, g_n\rangle\; d\lambda+\frac{1}{n}.  $$
Since $(\phi_n(t) f_n(t))$ is weakly null for each $t\in K$, $(\phi_n f_n)$ is weakly null in $C(K, X)$ by Lemma \ref{le1}.

Recall that $T^*(y^*)=m_{y^*}\in rcabv(\Sigma, X^*)$, $y^*\in Y^*$. By Theorem \ref{Din}, 
%$\langle T(f), y^* \rangle =\int_K \langle g_{y^*}, f \rangle \;d\lambda$,  $y^*\in Y^*$, $f\in C(K,X)$. 
%for every $X$-valued $\Sigma$ simple function $f$ defined  on $K$, and thus for every $f\in C(K,X)$.
$\langle T(f), y_n^* \rangle =\langle m_{y_n^*}, f \rangle =\int_K \langle g_n, f \rangle \;d\lambda$, for all $f\in C(K,X)$.
Then 
\begin{align*} 
\int_K |\langle f_n, g_n\rangle|\; d\lambda &\le \int_K   \langle \phi_n f_n, g_n \rangle \; d\lambda +\frac{1}{n} \\
  &= \langle T(\phi_n f_n), y_n^*\rangle +\frac{1}{n}  \to 0,
\end{align*}
since $T$ is weak Dunford-Pettis.

$b)\Rightarrow a)$  Let $(f_n)$ be a weakly null sequence in $C(K, X)$ and $(y_n^*)$ be a weakly null sequence in $Y^*$. Without loss of generality assume $(y_n^*)$ is in $B_{Y^*}$. Then
$$|\langle  T(f_n), y_n^* \rangle |=|\int_K \langle f_n, g_n\rangle\; d\lambda|  \le \int_K |\langle f_n, g_n\rangle| \; d\lambda \to 0, $$
 and thus $T$ is  weak Dunford-Pettis.

\end{proof}

%In the proof of the next lemma we will need the following result.   

We will need the following result. 

\begin{lemma} \label{wpsum} (\cite[Proposition 2.1]{SIA})
Let $1<p<\infty$ and $(x_n)$ be a sequence in Banach space $X$. Then the following are equivalent:

(a) The sequence $(x_n)$ is weakly $p$-summable in $X$.

(b) The series $\sum_{n=1}^\infty a_n x_n$ converges unconditionally for all $(a_n)\in \ell_{p'}$. 

(c) There exists an operator $T\in L(\ell_{p'}, X)$ such that $T(e_n)=x_n$, where $(e_n)$ is the unit vector basis of $\ell_{p'}$.
\end{lemma}

\begin{lemma}\label{le8}
(i) Let $1\le p<\infty$. If $(f_n)$ is a weakly $p$-summable sequence in $C(K, X)$ and $(\phi_n)$ is bounded in $C(K)$, then $(\phi_n f_n)$ is 
weakly $p$-summable sequence in $C(K, X)$. 

(ii) Let $1\le p\le \infty$. If $(x_n)$ is weakly $p$-summable in $X$ and $(\phi_n)$ is bounded in $C(K)$, then $(\phi_n x_n)$ is 
weakly $p$-summable sequence in $C(K, X)$. 
\end{lemma}

\begin{proof} Assume without loss of generality that $\|\phi_n\|_\infty\le 1$ for each $n$. 

(i) Let $1<p<\infty$. 
It is enough to prove that  $\sum_{n=1}^\infty a_n \phi_n f_n$ is unconditionally convergent, for all $(a_n)\in \ell_{p'}$ (by Lemma \ref{wpsum}).
 Let $(a_n)$ be a sequence in $B_{\ell_{p'}}$. Suppose that the subseries $\sum_{i=1}^\infty a_{n_i} \phi_{n_i} f_{n_i}$ is not convergent. Then there exists an $\epsilon>0$ and sequences $(p_n)$, $(q_n)$ of positive integers such that 
$p_1\le q_1<p_2\le q_2<p_3\le q_3<\cdots $, and $ \|\sum_{i=p_n}^{q_n} a_{n_i}\phi_{n_i} f_{n_i}\|_\infty >\epsilon$ for each $n$.
For each $n$, let $t_n\in K$ such that 
$$\|\sum_{i=p_n}^{q_n} a_{n_i} \phi_{n_i} f_{n_i} \|_\infty\le \|\sum_{i=p_n}^{q_n} a_{n_i}\phi_{n_i}(t_n) f_{n_i}(t_n)\|_\infty+\frac{\epsilon}{n}. $$
%\frac{\epsilon}{2^n}. 

Let $b_{n_i}=a_{n_i}\phi_{n_i}(t_n)$, for $p_n\le i\le q_n$, and $b_i=0$, for $i\not \in \cup_{n=1}^\infty \{p_n, \cdots q_n\}$. 
Then  
$$ \sum_{n=1}^\infty \sum_{i=p_n}^{q_n} |b_{n_i}|^{p'}\le \sum_{n=1}^\infty \sum_{i=p_n}^{q_n} |a_{n_i}|^{p'}\le \sum_{i=1}^\infty |a_i|^{p'}\le 1, $$
and thus $b=(b_i)\in B_{\ell_{p'}}$. We have 
\begin{align*}
 \|\sum_{i=p_n}^{q_n} a_{n_i}\phi_{n_i} f_{n_i}\|_\infty &\le \|\sum_{i=p_n}^{q_n} b_{n_i} f_{n_i}(t_n)\|_\infty +\frac{\epsilon}{n}\\
&\le  \|\sum_{i=p_n}^{q_n} b_{n_i} f_{n_i}\|_\infty +\frac{\epsilon}{n}\to 0, \qquad n\to \infty,
\end{align*}
because $\sum_{i=1}^\infty  b_{n_i} f_{n_i}  $ is  convergent. This contradiction concludes the proof.  

Let $p=1$. If $\sum_n f_n$ is a wuc series in $C(K, X)$  and $(\phi_n)$ is bounded in $C(K)$, then $\sum_n \phi_n f_n$ is 
wuc in $C(K, X)$. 

(ii) 
% Assume without loss of generality that $\|\phi_n\|_\infty\le 1$ for each $n$. 
Let $1<p<\infty$. It is enough to prove that  $\sum_{n=1}^\infty a_n \phi_n x_n$ is unconditionally convergent, for all $(a_n)\in \ell_{p'}$ (by Lemma \ref{wpsum}). Let $(a_n)$ be a sequence in $B_{\ell_{p'}}$. Suppose that the subseries $\sum_{i=1}^\infty a_{n_i} \phi_{n_i} x_{n_i}$ is not convergent. Then there exists an $\epsilon>0$ and sequences $(p_n)$, $(q_n)$ of positive integers such that 
$p_1\le q_1<p_2\le q_2<p_3\le q_3<\cdots $, and $ \|\sum_{i=p_n}^{q_n} a_{n_i}\phi_{n_i} x_{n_i} \|_\infty >\epsilon$ for each $n$.
For each $n$, let $t_n\in K$ such that 
$$\|\sum_{i=p_n}^{q_n}  a_{n_i}\phi_{n_i} x_{n_i} \|_\infty\le \|\sum_{i=p_n}^{q_n} a_{n_i}\phi_{n_i}(t_n) x_{n_i}\|+\frac{\epsilon}{n}. $$

%As in Lemma \ref{le8}, l
%As above, l
Let $b_{n_i}=a_{n_i} \phi_{n_i}(t_n)$, for $p_n\le i\le q_n$, and $b_i=0$, for $i\not \in \cup_{n=1}^\infty \{p_n, \cdots q_n\}$. 
Then  $b=(b_{i})\in B_{\ell_{p'}}$ and 
$$
 \|\sum_{i=p_n}^{q_n} a_{n_i}\phi_{n_i} x_{n_i} \|_\infty \le  \|\sum_{i=p_n}^{q_n} b_{n_i} x_{n_i} \| +\frac{\epsilon}{n}\to 0, \qquad n\to \infty,
$$
because $\sum_{i=1}^\infty  b_{n_i} x_{n_i}  $ is convergent. This contradiction concludes the proof.  

Let $p=1$. If $\sum_n x_n$ is a wuc series in $X$ and $(\phi_n)$ is bounded in $C(K)$, then $\sum_n \phi_n x_n$ is wuc in $C(K, X)$. Indeed, since $\sum_n x_n$ is a wuc, there exists an $M>0$ such that 
$$\|\sum_{i=1}^n \alpha_i x_i\|\le M \max\{|\alpha_i|: 1\le i\le n\}, $$
for each $n\in \mathbb{N}$ and $\alpha_1, \alpha_2, \ldots ,\alpha_n\in \mathbb{R}$ \cite[Proposition 4.3.9, p. 390]{Megginson}. For each $n\in \mathbb{N}$ and $t\in K$, 
$$\|\sum_{i=1}^n \phi_i(t) x_i\|\le M \max\{|\phi_i(t)|: 1\le i\le n\}\le M. $$
Then  
$$\|\sum_{i=1}^n \phi_i x_i\|_\infty \le M, $$
for each $n\in \mathbb{N}$, and  $\sum_n \phi_n x_n$ is wuc in $C(K, X)$.

Let $p=\infty$. If $(x_n)$ is weakly null in $X$ and $(\phi_n)$ is bounded in $C(K)$, then $(\phi_n x_n)$ is weakly null in $C(K, X)$ by Lemma \ref{le1}. 
%Let $p=1$. Suppose $\sum_n x_n$ is a wuc series in $C(K, X)$.  For each $t\in K$, $\sum \phi_i(t) x_i $ is wuc in $X$. Let $S_n=\sum_{i=1}^n \phi_i x_i$. Since $(S_n)$ is bounded in $C(K,X)$ and  $(S_n(t))$ is weakly Cauchy for each $t\in K$, $(S_n)$ is weakly Cauchy by Lemma \ref{le1}. Hence  $\sum_n \phi_n x_n$ is wuc in $C(K, X)$. 
\end{proof}

Lemma \ref{le8}, (ii) will be used later in the paper.

%The following result is similar to \cite[Theorem 1.5]{BRS}. 

\begin{theorem} \label{BomRS2}
Let $1\le p< \infty$. Let $m\leftrightarrow T:C(K, X)\to Y$ be a strongly bounded operator whose representing measure $m$ has a control measure $\lambda$. 
The following are equivalent:

a) $T$ is  weak $p$-convergent (resp. weak$^*$ $p$-convergent). 
%, weak $p$-convergent, weak$^*$ $p$-convergent

%
b) For every weakly $p$-summable sequence $(f_n)$ in $C(K, X)$ and every weakly null (resp. $w^*$-null) sequence $(y_n^*)$ in $B_{Y^*}$, we have 
$$\lim_{n\to \infty} \int_K |\langle f_n, g_n\rangle|\; d\lambda=0,  $$
where $g_n$ is the function corresponding to $T^*(y_n^*)$ by Theorem \ref{Din}. 
\end{theorem}

\begin{proof} 
The proof is similar to that of Theorem \ref{BomRS1}, using Lemma \ref{le8}. 
%that for every weakly $p$-summable sequence $(f_n)$ in $C(K, X)$ and every  bounded sequence $(\phi_n)$ in $C(K)$, the sequence $(\phi_n f_n)$ is also weakly $p$-summable in $C(K, X)$ (by Lemma \ref{le8}). 
\end{proof}

A positive Radon measure is discrete if every set of positive measure contains an atom (\cite{L}, Ch 2, Sect. 8). In this case, the measure is of the form $\sum_i a_i \delta_{t_i}$, with $\sum |a_i|<\infty$. In particular, it is concentrated on a countable set of its atoms.

%\vspace{0.1in}
%\noindent\textbf{ Observation.} If $(f_n)$ is  is weakly $p$-summable  in $C(K,X)$, then  $(f_n(t))$ is weakly  $p$-summable, for every $t\in K$.

%To see this, let $t\in K$. For each $x^*\in X^*$, let  $\delta_t(\cdot) x^* \in C(K,X)^*$,
 %$$ \langle \delta_t(\cdot) x^*, f\rangle=\int_K f d(\delta_t(\cdot) x^*)= \langle x^*, f(t)\rangle, \qquad f\in C(K, X).$$
%Since $(f_n)$  is weakly $p$-summable, 
 %$$\sum_{n=1}^\infty |\langle \delta_t(\cdot) x^*, f_n\rangle|^p =\sum_{n=1}^\infty | \langle x^*, f_n(t)\rangle|^p<\infty. $$
%Therefore  $(f_n(t))$ is weakly  $p$-summable. 

%$1<p<\infty$
\begin{theorem} \label{BomRS3}
Let $1\le p<\infty$. Let $m\leftrightarrow T:C(K, X)\to Y$ such that
% whose representing measure $m$ has a control measure $\lambda$. 
%The following are equivalent:

a) $m$ is  strongly bounded  and admits a discrete control measure $\lambda$.

b) For every   $A\in \Sigma$, $m(A):X\to Y$ is   weak Dunford-Pettis (resp. weak$^*$ Dunford-Pettis, weak $p$-convergent, weak$^*$ $p$-convergent).

Then $T$ is  weak Dunford-Pettis (resp.  weak$^*$ Dunford-Pettis, weak $p$-convergent, weak$^*$ $p$-convergent). 
\end{theorem}

%| |  | |

\begin{proof} We will prove the result for weak Dunford-Pettis operators. 
Let $(x_n)$ be a weakly null sequence in $X$  and $(y_n^*)$ be a weakly null sequence in $Y^*$. 
Without loss of generality assume $(y_n^*)$ is in $B_{Y^*}$. For each  $n\in \mathbb{N}$, let  $g_n$ be the function corresponding to $T^*(y_n^*)$ by Theorem \ref{Din}. 
 For every $A\in \Sigma$, $m(A):X\to Y$ is   weak Dunford-Pettis, and thus
$$ \lim_{n\to \infty} \int_A \langle x_n, g_n\rangle  \; d\lambda =\lim_{n\to \infty}  \langle m(A)(x_n), y_n^*\rangle = 0. $$
Then
 $$\lim_{n\to \infty} \, \langle x_n, g_n(t)\rangle =0, $$
for every $t\in K$ such that $\lambda(t)>0$.

Let $(f_n)$ be weakly null in $C(K, X)$. Then for every $t\in K$,  $(f_n(t))$ is weakly null in $X$ (by Lemma \ref{le1}). Therefore
$$\lim_{n\to \infty} \, \langle f_n(t), g_n(t)\rangle =0, $$
for every $t\in K$ such that $\lambda(t)>0$. 

The sequence $(f_n)$ is bounded and the set $\{|g_n|: n\in \mathbb{N}\}$ is uniformly integrable (since it is relatively weakly compact in $L^1(\lambda)$ \cite[p. 76]{DU}).  By Vitali's Theorem we obtain
$$ \lim_{n\to \infty}  \int_K |\langle f_n, g_n\rangle |\; d\lambda=0.$$
Therefore $T$ is  weak Dunford-Pettis by Theorem \ref{BomRS1}. 

The proof for weak $p$-convergent (resp. weak$^*$ $p$-convergent) operators is similar, using Theorem \ref{BomRS2} instead of Theorem \ref{BomRS1} and the fact that if $(f_n)$ is weakly $p$-summable in $C(K, X)$, then for every $t\in K$,  $(f_n(t))$ is weakly $p$-summable in $X$. 
%The proofs for the other cases are similar. 
\end{proof}

%$1<p< \infty$. 
\begin{corollary} \label{Cem1}
Let $1\le p< \infty$. Suppose that $K$ is a dispersed compact Hausdorff space $m\leftrightarrow T: C(K,X)\to Y$ is a strongly bounded operator. Then $T$  is  weak Dunford-Pettis (resp. weak$^*$ Dunford-Pettis, weak $p$-convergent, weak$^*$ $p$-convergent) if and only if for each $A\in \Sigma$, $m(A):X\to Y$ is   weak Dunford-Pettis (resp. weak$^*$ Dunford-Pettis, weak $p$-convergent, weak$^*$ $p$-convergent).
\end{corollary}

% Suppose that $K$ is a dispersed compact Hausdorff space $m\leftrightarrow T: C(K,X)\to Y$ is a strongly bounded operator such that for each $A\in \Sigma$, $m(A):X\to Y$ is   weak Dunford-Pettis (resp. weak$^*$ Dunford-Pettis, weak $p$-convergent, weak$^*$ $p$-convergent). Then  $T$  is  weak Dunford-Pettis (resp. weak$^*$ Dunford-Pettis, weak $p$-convergent, weak$^*$ $p$-convergent). 

\begin{proof}
Suppose that for each $A\in \Sigma$, $m(A):X\to Y$ is weak Dunford-Pettis. Then $T$ is weak Dunford-Pettis by Theorem \ref{BomRS3}, since every Radon measure on a compact dispersed space $K$ is discrete \cite[Ch. 2, Sect. 8]{L}.

If $T$  is  weak Dunford-Pettis, then for each $A\in \Sigma$, $m(A):X\to Y$ is   weak Dunford-Pettis by Corollary \ref{mwDP}. 
\end{proof}

 A Banach space $X$ has the $DP^*$-property ($DP^*P$) if all weakly compact sets in $X$ are limited  \cite{CGL}.

% \cite{IGCC}
 A Banach space $X$ has the $DP^*P$ if and only if  $x_n^*(x_n)\to 0$ for all weakly null sequences $(x_n)$ in $X$ and $w^*$-null sequences $(x_n^*)$ in $X^*$. A  space  $X$ has the $DP^*P$ if and only if every operator $T:X\to c_0$ is completely continuous. 
\cite{CGL}. 
If $X$ has the $DP^*P$, then it has the $DPP$. 
If $X$ is a Schur space or if $X$ has the $DPP$ and the Grothendieck property (weak and weak$^*$ convergence of sequences in $X^*$ coincide), then $X$ has the $DP^*P$.

Let $1\le p\le \infty$. A Banach space $X$ has the \emph{Dunford-Pettis property of order p}  $(DPP_p)$ if every weakly compact operator $T:X \to Y$ is $p$-convergent, for any Banach space $Y$ \cite{CS}. Equivalently, $X$ has the $DPP_p$ if and only if $x_n^*(x_n)\to 0$ whenever $(x_n)$ is  weakly $p$-summable  in $X$ and $(x_n^*)$ is  weakly null  in $X^*$ \cite[Proposition 3.2]{CS}.

%If $X$ has the $DPP_p$, then it has the $DPP_q$, if $q<p$. 
The $DPP_\infty$ is precisely the $DPP$, and every Banach space  has the $DPP_1$. $C(K)$ spaces and $L_1$ have the $DPP$, and thus the $DPP_p$ for all $p$.

Let $1\le p\le \infty$. A Banach space $X$ has the $DP^*$-\emph{property of order} $p$  ($DP^*P_p$) if all weakly $p$-compact sets in $X$ are limited \cite{FZ2}.
%\cite{FZ}. 
%A  space  $X$ has the $DP^*P_p$ if and only if $L(X, c_0)=C_p(X, c_0)$ \cite{FZ}.
The Banach space $X$ has the $DP^*P_p$ if and only if $x_n^*(x_n)\to 0$ for all weakly $p$-summable sequences $(x_n)$ in $X$ and $w^*$-null sequences $(x_n^*)$ in $X^*$ \cite[Theorem 2.7]{FZ2}.

%If $X$ has the $DP^*P_q$, then it has the $DP^*P_p$, if $q>p$. 
The $DP^*P_\infty$ is precisely the $DP^*P$ and every Banach space has the $DP^*P_1$. If $X$ has the  $DP^*P$, then  $X$ has the $DP^*P_p$, for all $1 \le p \le \infty$. 
If $X$ has the $DP^*P_p$, then $X$ has the $DPP_p$.

%Let $X$ and $Y$ be Banach spaces.

We note that if $X$ has the $DPP$ (resp. $DP^*P$,  $DPP_p$, $DP^*P_p$), then every operator $T:X\to Y$ is weak Dunford-Pettis (resp. weak$^*$ Dunford-Pettis, weak $p$-convergent, weak$^*$ $p$-convergent), for $1<p<\infty$ \cite{GRen}, \cite{FZ2}.
%if and only if the identity operator $i:X\to X$ is weak Dunford-Pettis (resp. weak$^*$ Dunford-Pettis, weak $p$-convergent, weak$^*$ $p$-convergent) \cite{GRen}, \cite{FZ2}.  

%It is known that if $K$ is a dispersed compact Hausdorff space, then $C(K, X)$ has the $DPP$ (resp. the $DPP_p$) if and only if $X$ has the $DPP$ (resp. the $DPP_p$) \cite{BCCam}, \cite{CS}. It follows that if $K$ is dispersed and $X$ has the $DPP$ (resp. the $DPP_p$), then every operator $T:C(K, X)\to Y$ is weak Dunford-Pettis (resp. weak $p$-convergent). 

\begin{corollary}
 Let $1<p<\infty$ and let $K$ be a dispersed compact Hausdorff space.

(i) If $X$ has the $DPP$ (resp. $DPP_p$), then every operator $T:C(K, X)\to Y$ is weak Dunford-Pettis (resp. weak $p$-convergent), for every Banach space $Y$. 

(ii) If $X$ has the $DP^*P$ (resp. $DP^*P_p$), then every strongly bounded  operator $T: C(K,X)\to Y$ is weak$^*$ Dunford-Pettis (resp. weak$^*$ $p$-convergent), for every Banach space $Y$. 

(iii) If $X$ has the $DP^*P$ (resp. $DP^*P_p$), then every strongly bounded  operator $T: C(K,X)\to c_0$ is completely continuous (resp. $p$-convergent). 
\end{corollary}

\begin{proof}

(i) If $K$ is dispersed and $X$ has the $DPP$ (resp. $DPP_p$), then $C(K, X)$ has the $DPP$ (resp. $DPP_p$) \cite{BCCam}, \cite{CS}. 
Then every operator $T:C(K, X)\to Y$ is weak Dunford-Pettis (resp. weak $p$-convergent).

(ii) Let $m\leftrightarrow T: C(K,X)\to Y$ be a strongly bounded  operator. Since   $X$ has the $DP^*P$ (resp. $DP^*P_p$), $m(A):X \to Y$ is weak$^*$ Dunford-Pettis (resp.  weak$^*$ $p$-convergent), for every $A\in \Sigma$.
%  \cite{GRen}, \cite{FZ2}. 
Therefore $T$ is weak$^*$ Dunford-Pettis (resp.  weak$^*$ $p$-convergent), by Corollary \ref{Cem1}. 

(iii) Let $T: C(K,X)\to c_0$ be a strongly bounded operator. By (ii), $T$ is weak$^*$ Dunford-Pettis (resp.  weak$^*$ $p$-convergent). Since $c_0$ is separable, $T$ is completely continuous (resp. $p$-convergent) by \cite[Corollary 10]{GRen} (resp. \cite[Corollary 23]{GRen}). 
\end{proof}

If $T:C(K, X)\to Y$ is an operator such that $T^*$ is weakly precompact, then $T$ is unconditionally converging  \cite[Theorem 1]{BLNach}, hence strongly bounded \cite{ID}. 

If $K$ is a compact Hausdorff space and $\phi\in C(K)$, for every Banach space $X$ the linear map $S_{\phi}:C(K, X)\to C(K, X)$ defined by 
$S_{\phi}(f)=\phi f$ is continuous. For a measure $m\in rcabv(\Sigma, X^*)=C(K,X)^*$, let  $S_{\phi}^*(m)=\phi m$ \cite{FBom}. 

%$1\le p\le \infty$
\begin{lemma} \label{le17}
Let $1\le  p< \infty$. If $(\phi_n)$ is bounded in $C(K)$ and $(f_n)$ is DP weakly null (resp. DP weakly $p$-summable)  in $C(K, X)$, then $(\phi_n f_n)$ is DP weakly null (resp. DP weakly $p$-summable)  in $C(K, X)$.
\end{lemma}

\begin{proof}
Suppose $(f_n)$ is DP weakly null in $C(K, X)$. 
%Assume without loss of generality that $\|\phi_n\|_\infty\le 1$ for each $n$. 
For each $t\in K$, $(\phi_n(t) f_n(t))$ is weakly null in $X$, and thus $(\phi_n f_n)$ is  weakly null in $C(K, X)$ by Lemma \ref{le1}. 
  Let $T:C(K, X)\to Y$ be an operator so that $T^*$ is weakly precompact. 
	%We show $\|T(\phi_n f_n)\| \to 0$. 
	%It is enough to show that $\|T(\phi_n f_n)\| \to 0$,  \cite[Corollary 4]{IGMon}. 
	For each $n$, let $y_n^*\in Y^*$ such that $\|y_n^*\|= 1$ and $\langle T(\phi_n f_n), y_n^*\rangle =\|T(\phi_n f_n)\| $. 
Let $m_n=T^*(y_n^*)$ for each $n$. Since $T^*$ is weakly precompact, the set $\{ m_n: n \in \mathbb{N} \}$ is weakly precompact. Hence  $\{\phi_n m_n: n \in \mathbb{N} \}$ is weakly precompact \cite[Lemma 3]{FBom}. By passing to a subsequence, we can assume that $(\phi_n m_n)$ is weakly Cauchy. Since $(f_n)$ is DP weakly null, 
$$ \langle \phi_n f_n, m_n \rangle= \langle  f_n, \phi_n m_n \rangle\to 0, $$
by \cite[Proposition 22]{GPQM}.
Therefore $\|T(\phi_n f_n)\| \to 0$, and thus  $(\phi_n f_n)$ is   DP  in $C(K, X)$ by  \cite[Corollary 4]{IGMon}. 

Now suppose $(f_n)$ is DP weakly $p$-summable in $C(K, X)$. Then $(\phi_n f_n)$ is weakly $p$-summable in $C(K, X)$ by Lemma  \ref{le8}
 and $(\phi_n f_n)$ is DP by the previous argument.
\end{proof}

 An operator $T:X\to Y$ is \emph{pseudo weakly compact (pwc)} (or \emph{Dunford-Pettis completely continuous (DPcc)}) if it takes DP weakly null  sequences in $X$ into norm null sequences in $Y$ \cite{PVWY}, \cite{WC}. 

Let $1\le p< \infty$. An operator $T:X\to Y$ is called \emph{DP $p$-convergent} if it takes DP weakly $p$-summable sequences to norm null sequences \cite{GLBCl}.

An operator  $T:X\to Y$ is called \emph{limited completely continuous (lcc)} if  $T$ maps limited weakly null  sequences to norm null sequences \cite{SM}.

Let $1\le p< \infty$. An operator $T:X\to Y$ is called \emph{limited $p$-convergent} if it takes limited weakly $p$-summable sequences in $X$ to norm null ones in $Y$ \cite{FZ2}. 

%The case  $p=\infty$ could have been included in the previous definitions if we had considered that the limited (resp. DP) $\infty$-convergent operators are precisely the limited completely continuous (resp. pseudo weakly compact) operators.  

%If $m\leftrightarrow T: C(K,X)\to Y$ is a  pseudo weakly compact (resp. DP $p$-convergent operator), then  $T$ is unconditionally converging by \cite[Proposition 14]{PVWY} (resp. \cite[ Proposition 14]{GPQM}), hence strongly bounded \cite{ID}. Further, $m(A):X\to Y$ is pseudo weakly compact (resp. DP $p$-convergent) for each $A\in \Sigma$ \cite{GPQM}.

If $m\leftrightarrow T: C(K,X)\to Y$ is a  pseudo weakly compact (resp. DP $p$-convergent operator), then $T$ is strongly bounded and $m(A):X\to Y$ is pseudo weakly compact (resp. DP $p$-convergent) for each $A\in \Sigma$ \cite{GPQM}. 
%\cite[Corollary 17]{GPQM}

\begin{theorem} \label{BomRS1-pwc}
Let $1\le p<\infty$. Let $m\leftrightarrow T:C(K, X)\to Y$ be an operator whose representing measure $m$ has a control measure $\lambda$. 
The following are equivalent:

a) $T$ is pseudo weakly compact (resp.  DP $p$-convergent).

b) For every  DP weakly null (resp. DP weakly $p$-summable)  sequence $(f_n)$ in $C(K, X)$ and every sequence $(y_n^*)$ in $B_{Y^*}$, we have 
$$\lim_{n\to \infty} \int_K |\langle f_n, g_n\rangle|\; d\lambda=0,  $$
where $g_n$ is the function corresponding to $T^*(y_n^*)$ by Theorem \ref{Din}. 
\end{theorem}

\begin{proof} 
%The proof is similar to the proof of \cite[Theorem 1.4]{BRS}, using  Lemma \ref{le17}.
%the fact that if $(\phi_n)$ is bounded in $C(K)$ and $(f_n)$ is DP weakly null (resp. DP weakly $p$-summable)  in $C(K, X)$, then $(\phi_n x_n)$ is DP weakly null (resp. DP weakly $p$-summable)  in $C(K, X)$ by Lemma \ref{le17}.
%  We only prove the result for  pseudo weakly compact operators; the other case is similar. 
$a)\Rightarrow b)$ Suppose $T$ is pseudo weakly compact.  Let $(f_n)$ be  a  DP weakly null sequence in $C(K, X)$ and $(y_n^*)$ be a sequence in $B_{Y^*}$. Let $\phi_n$ be a scalar continuous function on $K$ such that $\|\phi_n\|\le 1 $ and 
$$ \int_K |\langle f_n, g_n\rangle|\; d\lambda \le  \int_K \phi_n \langle f_n, g_n\rangle\; d\lambda+\frac{1}{n}.  $$
Then  $(\phi_n f_n)$ is DP weakly null in $C(K, X)$ by Lemma \ref{le17} and
\begin{align*} 
\int_K |\langle f_n, g_n\rangle|\; d\lambda &\le \langle T(\phi_n f_n), y_n^*\rangle +\frac{1}{n} \le \|T(\phi_n f_n)\|+\frac{1}{n} \to 0,
\end{align*}
since $T$ is pseudo weakly compact.

$b)\Rightarrow a)$  Let $(f_n)$ be a DP weakly null sequence in $C(K, X)$ and let $(y_n^*)$ be a sequence in $B_{Y^*}$ such that
$\langle y_n^*, T(f_n) \rangle =\|T(f_n)\|$ for each $n$. Then
$$\|T(f_n)\|=
\int_K \langle f_n, g_n\rangle\; d\lambda  \le \int_K |\langle f_n, g_n\rangle| \; d\lambda \to 0, $$
 and thus $T$ is pseudo weakly compact.
\end{proof}

%With the notations of Theorem \ref{BomRS1-pwc}, 
\noindent\textbf{Remark.} Let $m\leftrightarrow T:C(K, X)\to Y$ be a strongly bounded operator whose representing measure $m$ has a control measure $\lambda$. Suppose that for every  limited weakly null (resp. limited weakly $p$-summable)  sequence $(f_n)$ in $C(K, X)$ and every sequence $(y_n^*)$ in $B_{Y^*}$, we have 
$$\lim_{n\to \infty} \int_K |\langle f_n, g_n\rangle|\; d\lambda=0,  $$
where $g_n$ is the function corresponding to $T^*(y_n^*)$ by Theorem \ref{Din}. Then $T$ is limited completely continuous (resp.  limited $p$-convergent). The proof is similar to the proof of $b)\Rightarrow a)$ in Theorem \ref{BomRS1-pwc}. 

\begin{theorem} \label{BomRS3-pwc}
Let $1\le p<\infty$. Let $m\leftrightarrow T:C(K, X)\to Y$ such that

a) $m$ is  strongly bounded  and admits a discrete control measure $\lambda$.

b) For every $A\in \Sigma$, $m(A):X\to Y$ is pseudo weakly compact (resp.  DP $p$-convergent).

Then $T$ is  pseudo weakly compact (resp.  DP $p$-convergent). 
\end{theorem}

\begin{proof} We only prove the result for  pseudo weakly compact operators; the other case is similar. 
Let $(x_n)$ be a DP weakly null sequence in $X$  and $(y_n^*)$ be a  sequence in $B_{Y^*}$. 
%Without loss of generality assume $(y_n^*)$ is in $B_{Y^*}$.
 For every $A\in \Sigma$, $m(A):X\to Y$ is pseudo weakly compact, and thus
$$ | \int_A \langle x_n, g_n\rangle  \; d\lambda |=| \langle m(A)(x_n), y_n^*\rangle|\le \|m(A)(x_n)\|\to 0. $$
Then
 $$\lim_{n\to \infty} \, \langle x_n, g_n(t)\rangle =0, $$
for every $t\in K$ such that $\lambda(t)>0$. 
%since $\lambda$ is discrete

Let $(f_n)$ be a DP weakly null sequence in $C(K, X)$. Then for every $t\in K$,  $(f_n(t))$ is DP weakly null in $X$ (by the proof of  \cite[Theorem 20]{GPQM}). Hence
$$\lim_{n\to \infty} \, \langle f_n(t), g_n(t)\rangle =0, $$
for every $t\in K$ such that $\lambda(t)>0$. 

%The sequence $(f_n)$ is bounded and the set $\{|g_n|: n\in \mathbb{N}\}$ is uniformly integrable (since it is relatively weakly compact in $L^1(\lambda)$ \cite[p. 76]{DU}).  By Vitali's Theorem we obtain
As in the proof of Theorem \ref{BomRS3}, by Vitali's Theorem we obtain
$$ \lim_{n\to \infty}  \int_K |\langle f_n, g_n\rangle |\; d\lambda=0.$$
Therefore $T$ is pseudo weakly compact by Theorem \ref{BomRS1-pwc}. 
\end{proof}

\noindent
\textbf{Remark.} Every Radon  measure on a compact dispersed space $K$ is discrete \cite[Ch. 2, Sect. 8]{L}. Therefore Theorem \ref{BomRS3-pwc} extends \cite[Theorem 20]{GPQM}. 

\begin{lemma}\label{lemmalcc}
Let $1\le p<\infty$. 
If $(f_n)$ is a limited weakly null (resp. limited weakly $p$-summable) sequence in $C(K, X)$, then  $(f_n(t))$ is  limited weakly null (resp. limited weakly $p$-summable) for each $t\in K$. 
\end{lemma}

\begin{proof}
Let $(f_n)$ be a limited weakly null sequence in $C(K, X)$. Then $(f_n(t))$ is a weakly null sequence in $X$  for each $t\in K$ (by Lemma \ref{le1}). Let $t\in K$ . Consider the linear map $T:X^*\to C(K, X)^*$, $T(x^*)=\delta_t(\cdot) x^*$. 
%For each $x^*\in X^*$, $\delta_t(\cdot) x^* \in C(K,X)^*$,
 $$ \langle T(x^*), f\rangle=\int_K f d(\delta_t(\cdot) x^*)= \langle x^*, f(t)\rangle, \qquad f\in C(K,X).$$
% \langle \delta_t(\cdot) x^*, f\rangle
If $f\in C(K, X)$, then 
$\langle T^*(f), x^* \rangle= \langle x^*, f(t)\rangle$, and $T^*(f)=f(t)$. Hence $T^*$ maps $C(K, X)$ into $X$, and thus
$T$ takes $w^*$-null sequences to $w^*$-null sequences. 
%the sequence
If $(x_n^*)$ is a $w^*$-null sequence in $X^*$, then  $(\delta_t(\cdot)x_n^*)$ is $w^*$-null in $C(K, X)^*$, and 
$$ \langle \delta_t(\cdot) x_n^*, f_n\rangle=\int_K f_n d(\delta_t(\cdot) x_n^*)= \langle x_n^*, f_n(t)\rangle \to 0.$$
Hence $(f_n(t))$ is  limited   for each $t\in K$.
%in $X$
\end{proof}

\begin{theorem} \label{BomRS3-lcc}
Let $1\le p<\infty$. Let $m\leftrightarrow T:C(K, X)\to Y$ such that

a) $m$ is  strongly bounded  and admits a discrete control measure $\lambda$.

b) For every $A\in \Sigma$, $m(A):X\to Y$ is limited completely continuous (resp.  limited $p$-convergent).

Then $T$ is  limited completely continuous (resp.  limited $p$-convergent).
\end{theorem}

\begin{proof}
The proof is similar to the proof of Theorem \ref{BomRS3-pwc}, using Lemma \ref{lemmalcc} and the remark after Theorem \ref{BomRS1-pwc}. 
\end{proof}

%\vspace{0.1in}
%An operator  $T:X\to Y$ is called \emph{limited completely continuous (lcc)} if  $T$ maps limited weakly null  sequences to norm null sequences \cite{SM}.

%Let $1\le p< \infty$. An operator $T:X\to Y$ is called \emph{limited $p$-convergent} if it carries limited weakly $p$-summable sequences in $X$ to norm null ones in $Y$ \cite{FZ2}.

\begin{corollary} \label{Tlcc} Let $1\le p< \infty$.  Let $K$ be a dispersed compact Hausdorff space. 
Suppose  $m\leftrightarrow T: C(K,X)\to Y$ is a strongly bounded operator such that  for each $A\in \Sigma$, $m(A):X\to Y$ is limited completely continuous (resp.  limited $p$-convergent). Then $T$  is limited completely continuous (resp.  limited $p$-convergent).
\end{corollary}

\begin{proof} Theorem \ref{BomRS3-lcc} implies that 
$T$  is limited completely continuous (resp.  limited $p$-convergent), since every Radon  measure on a compact dispersed space $K$ is discrete \cite[Ch. 2, Sect. 8]{L}.
% We only prove the result for limited completely continuous operators. The other case is similar.
%Let $(f_n)$ be a limited weakly null sequence in $C(K, X)$. Then $(f_n(t))$ is a weakly null sequence in $X$  for each $t\in K$ (by Lemma \ref{le1}). Let $t\in K$ . Consider the linear map $T:X^*\to C(K, X)^*$, $T(x^*)=\delta_t(\cdot) x^*$. 
% $$ \langle T(x^*), f\rangle=\int_K f d(\delta_t(\cdot) x^*)= \langle x^*, f(t)\rangle, \qquad f\in C(K,X).$$
%If $f\in C(K, X)$, then $\langle T^*(f), x^* \rangle= \langle x^*, f(t)\rangle$, and $T^*(f)=f(t)$. Hence $T^*$ maps $C(K, X)$ into $X$, and thus $T$ takes $w^*$-null sequences to $w^*$-null sequences. If $(x_n^*)$ is a $w^*$-null sequence in $X^*$, then  $(\delta_t(\cdot)x_n^*)$ is $w^*$-null in $C(K, X)^*$, and 
%$$ \langle \delta_t(\cdot) x_n^*, f_n\rangle=\int_K f_n d(\delta_t(\cdot) x_n^*)= \langle x_n^*, f_n(t)\rangle \to 0.$$
%Hence $(f_n(t))$ is a  limited sequence in $X$  for each $t\in K$.
% For each $A\in \Sigma$,  $m(A)(\{f_n(t) : n\in \mathbb{N}\}$ is relatively compact (since $m(A):X\to Y$ is  limited completely continuous). Then $\{T(f_n): n\in \mathbb{N}\}$ is relatively compact by \cite[ Remark 1]{GLBCl}. Thus $\|T(f_n)\|\to 0$ and $T$ is limited completely continuous.
\end{proof}

%\vspace{0.1in}
In \cite{ABBL} it was given an alternative proof of the fact that if $m\leftrightarrow T: C(K,X)\to Y$ is a Dieudonn\'e (or weakly completely continuous) operator, then $m(A):X\to Y$  is a Dieudonn\'e   operator for each $A\in \Sigma$. We use a similar argument for  limited completely continuous (resp. limited $p$-convergent) operators. 

\begin{proposition} \label{closed}
Let $1\le p< \infty$. 
If $T_n:X\to Y$ is a limited completely continuous (resp. limited $p$-convergent, pseudo weakly compact, DP $p$-convergent) operator for each $n\in \mathbb{N}$ and $T:X\to Y$ is an operator such that $\|T_n-T\|\to 0$, then $T$ is limited completely continuous (resp. limited $p$-convergent, pseudo weakly compact, DP $p$-convergent). 
\end{proposition}

\begin{proof} 
Suppose $(x_n)$ is a limited weakly null sequence in $X$ and let $\epsilon>0$. Assume without loss of generality that $\|x_n\|\le 1$ for each $n$. 
%Since $\|T_k-T\|\to 0$, 
Let $k_0\in \mathbb{N}$ such that  $\|T_k-T\|<\epsilon $, $k\ge k_0$. Let $k\ge k_0$. Since $T_k$ is limited completely continuous, $\|T_k(x_n)\|\to 0$  as $n\to \infty$. Choose $n$ sufficiently large such that $\|T_k(x_n)\|<\epsilon$. Then
$$ \|T(x_n)\|\le \|(T_k-T)(x_n)\|+\|T_k(x_n)\| <2\epsilon, $$
 %\epsilon + \|T_k(x_n)\|
and thus $T$ is limited completely continuous. 
\end{proof}

\begin{theorem} \label{T2lcc}
Let $1\le p< \infty$. 
% Let $K$ be a compact Hausdorff space. 
Suppose  $m\leftrightarrow T: C(K,X)\to Y$ is a strongly bounded operator. If $T$  is limited completely continuous (resp.  limited $p$-convergent), then  $m(A):X\to Y$ is limited completely continuous (resp.  limited $p$-convergent) for each $A\in \Sigma$. 
\end{theorem}

\begin{proof} 
Suppose $m\leftrightarrow T: C(K,X)\to Y$ is a strongly bounded limited completely continuous operator  and $f\in C(K)$. Define $T_f:X\to Y$ by $$T_f(x)=T(fx), \qquad x\in X.$$
 Let $(x_k)$ be a limited weakly null sequence in $X$. Consider the map $S_f:X\to C(K, X)$, $S_f(x)=fx$, $x\in X$. Since $S_f$ takes limited weakly null sequences to limited weakly null sequences, $(f x_k)$ is limited weakly null. Therefore  
$\|T_f(x_k)\|=\|T(fx_k)\|\to 0$, and thus $T_f$ is limited completely continuous. 

Let $A\in \Sigma$. Choose a decreasing sequence of open sets $(O_n)$ and an increasing sequence of compact sets $(K_n)$ such that 
$K_n\subset A\subset O_n$ for each $n$ and $\tilde{m} (O_n\setminus K_n)\to 0$. Let $(f_n)$ be a sequence in $C(K)$ so that $\|f_n\|=1$, $\text{support}(f_n)\subset O_n$, and $f_n=1$ on $K_n$ for each $n$. 
%By Urysohn's lemma,
Define $T_n:X\to Y$, $T_n(x)=T(f_nx)$, $x\in X$, for each $n$. Then $T_n=T_{f_n}$ is limited completely continuous for each $n$. Let $T_A:X\to Y$ be defined by 
$T_A(x)=T^{**}(\chi_A x)=m(A)(x)$, $x\in X$. If $x\in B_X$, then 
\begin{align*}
 \|T_n(x)-T_A(x)\|  &=\|T^{**}(f_nx)-T^{**}(\chi_A x) \|\\
&= \|\int_{O_n\setminus K_n} (f_n x-\chi_A x) \, dm \|\le  \tilde{m} (O_n\setminus K_n)\to 0. 
\end{align*}
%$$
%\|T_n(x)-T_A(x)\|  =\|T^{**}(f_nx)-T^{**}(\chi_A x) \|\le  \tilde{m} (O_n\setminus K_n)\to 0. 
%$$

%\|T(f_nx)-T^{**}(\chi_A x) \|

It follows that $\|T_n-T_A\|\to 0$, and  $m(A)$ is limited completely continuous by Proposition \ref{closed}. 
\end{proof}

A subset $A$ of  $X^*$ is called an $L$-\emph{limited set}  \cite{SMA}  if  each limited weakly null sequence $(x_n)$ in $X$ tends to $0$ uniformly on $A$.

%A Banach space $X$ has the \emph{$L$-limited property} \cite{SMA}  if every $L$-limited subset of $X^*$ is relatively weakly compact.

A Banach space $X$ has the \emph{$L$-limited property} \cite{SMA} (resp. \emph{$wL$-limited property} \cite{IGExt}) if every $L$-limited subset of $X^*$ is relatively weakly compact  (resp. weakly precompact).

Let $1\le p<\infty$. A  subset $A$ of a dual space $X^*$ is called a  $p$-$L$-\emph{limited set} \cite{IGAOT}  if for every   limited weakly $p$-summable sequence $(x_n)$ in $X$,  $\sup_{x^* \in A } |x^* (x_n) | \to 0$. 
%or $L_p$-limited set \cite{DM2}) 

Let $1\le p< \infty$. A Banach space $X$ has the   \emph{$p$-$L$-limited} \cite{IGAOT} (resp.  \emph{$p$-$wL$-limited}) \emph{property} if every $p$-$L$-limited subset of $X^*$ is relatively weakly compact (resp. weakly precompact).

\begin{corollary} Let $K$ be a dispersed compact Hausdorff space. 
  
(i) If $X$ has the  $L$-limited (resp. $wL$-limited property)  property, then every strongly bounded limited completely continuous operator $T:C(K, X)\to Y$ is weakly compact (resp. weakly precompact).

%strongly bounded
(ii) Let $1\le p< \infty$. If $X$ has the $p$-$L$-limited (resp.  $p$-$wL$-limited) property, then every  strongly bounded limited $p$-convergent operator $T:C(K, X)\to Y$ is  weakly compact (resp. weakly precompact).

\end{corollary}

\begin{proof}
   
 (i) Let $m\leftrightarrow T: C(K,X)\to Y$ be a  strongly bounded limited completely continuous operator.  Then for each $A\in \Sigma$, $m(A): X\to Y$ is limited completely continuous by Theorem \ref{T2lcc}. Since $X$ has the $L$-limited  (resp. $wL$-limited) property,  $m(A):X\to Y$ has a  weakly compact (resp. weakly precompact) adjoint for each $A\in \Sigma$ \cite[Theorem 2.8]{SMA} (resp. \cite[Theorem 2]{IGExt}). Thus $m(A)$ is weakly compact  (resp. weakly precompact \cite[Corollary 2]{BLNach}). Hence $T$ is weakly compact (resp. weakly precompact by \cite[Theorem 7]{BCCam} (resp. \cite[Theorem 11]{GLBCl}).  

%strongly bounded 
(ii) Suppose $m\leftrightarrow T: C(K,X)\to Y$ is a  strongly bounded limited $p$-convergent operator. Then  $m(A): X\to Y$ is  limited $p$-convergent    for each $A\in \Sigma$ by Theorem \ref{T2lcc}. Since $X$ has the $p$-$L$-limited  (resp. $p$-$wL$-limited) property,    $m(A):X\to Y$ has a  weakly compact (resp. weakly precompact) adjoint for each $A\in \Sigma$ \cite[Theorem 3.10]{IGAOT}. Continue as above.
%Then $m(A)$ is weakly compact (resp. weakly precompact) for each $A\in \Sigma$. Hence $T$ is weakly compact \cite[Theorem 7]{BCCam} (resp. weakly precompact  \cite[Theorem 11]{GLBCl}).  
 %Suppose $X$ has the $p$-$(SR)$ property. By (ii), every  DP $p$-convergent operator $T: C(K,X)\to Y$ is weakly compact. Then $C(K, X)$ has the $p$-$(SR)$ property \cite[Theorem 3.10]{IGAOT}. 
\end{proof}

%The authors do not know whether $m\leftrightarrow T: C(K,X)\to Y$ is limited completely continuous (resp.  limited $p$-convergent) implies that $T$  is strongly bounded and for each $A\in \Sigma$, $m(A):X\to Y$ is limited completely continuous (resp.  limited $p$-convergent). 

% $1< p<\infty$

It is known that the following are equivalent:

(a) $K$ is dispersed. 

(b) There is a Banach space $X$ containing an isomorphic copy of $c_0$ (resp. $X$ does not have the Schur property) such that an operator $m\leftrightarrow T:C(K,X) \to c_0$ is unconditionally convergent (resp. completely continuous) if and only if $m$ is strongly bounded and  $m(A):X\to c_0$ is unconditionally convergent (resp. completely continuous)  for every $A\in \Sigma$  \cite[Theorem 9]{BCCam} (resp. \cite[Theorem 11]{BCCam}). We obtain a similar result for $p$-convergent operators. 

%$1< p<\infty$
\begin{theorem} \label{Th14}
Let $1\le p\le\infty$ and let $K$ be a compact Hausdorff space. The following statements are equivalent:

(i) $K$ is dispersed. 

%a strongly bounded  operator

(ii) For any pair of Banach spaces $X$ and $Y$, an operator $m\leftrightarrow T:C(K,X) \to Y$ is $p$-convergent if and only if $m$ is strongly bounded and $m(A):X\to Y$ is $p$-convergent for every $A\in \Sigma$.

%a strongly bounded
(iii) There is a Banach space $X$ with  $X\not \in C_p$ such that an operator $m\leftrightarrow T:C(K,X) \to c_0$ is  $p$-convergent if and only if $m$ is strongly bounded and $m(A):X\to c_0$ is   $p$-convergent  for every $A\in \Sigma$. 
\end{theorem}

\begin{proof}  The case $p=1$ is \cite[Theorem 9]{BCCam} (note that $X\in C_1$ if and only if $X$ does not contain a copy of $c_0$ by the Bessaga and Pe{\l}czy\'nski theorem \cite[Theorem 8, p. 45]{JDSS}). The case $p=\infty$ is \cite[Theorem 11]{BCCam}. 

Let $1<p<\infty$. $(i)  \Rightarrow (ii)$ by  \cite[Proposition 2.1]{CS} and \cite[Corollary 2.3]{CS}. 
$(ii)  \Rightarrow (iii)$ is clear.

$(iii)  \Rightarrow (i)$ Suppose that (iii) holds and $K$ is not dispersed.
Then there is a purely nonatomic regular probability Borel measure $\lambda$ on $K$ (\cite[Theorem 2.8.10]{L}). Now we can construct a Haar system $\{ A_i^n : \, 1\le i \le 2^n , n\ge 0\}$ in $\Sigma$ (that is, $A_1^0=K$, for each $n\in \mathbb{N}, \{ A_i^n: \, 1\le i\le 2^n \}$ is a partition of $K$, and $A_i^n=A_{2i-1}^{n+1} \cup A_{2i}^{n+1}$, $1\le i\le 2^n$, $n\ge 0$) such that $\lambda(A_i^n)=2^{-n}$, for $1\le i\le 2^n$, $n\ge 0$.
Let $X$ be a Banach space as in the statement of (iii). Since $X\not \in C_p$, 
 there exists  a  weakly $p$-summable  sequence $(x_n)$  in $X$ with $\|x_n \|=1$, $n \ge 0$. For each $n\ge 0$, choose $x_n^*\in X^*$ such that $x_n^*(x_n)=1= \|x_n^*\| $, and let $r_n=\sum_{i=1}^{2^n} (-1)^i \chi_{A_i^n}$. Then $(r_n)$ is orthonormal in $L^2(\lambda)$, and thus is weakly null in $L^1(\lambda)$. 
Define $T:C(K,X) \to c_0$ by
$$ T(f)= \left(\int_K \langle x_n^*,  f(t) \rangle \, r_n(t) \, d\lambda \right)_{n\ge 0}\,\, , \, f\in C(K,X). $$

We note that  $T(f)\in c_0$ for each $f\in C(K,X)$ and that  the representing measure $m$ of $T$ is given by
$$m(A)(x)= ( \langle x_n^*, x \rangle \int_A r_n(t) \, d\lambda)_{n\ge 0},$$
for $A\in \Sigma$, $x\in X$. 
Since $(\int_A r_n(t) \, d\lambda) \to 0$ for all $A\in \Sigma$ and $\{ \langle x_n^*, x \rangle: n\in \mathbb{N}, x\in B_X \}$ is bounded, it follows that  $m(A)(x)\in c_0$.  Further, $\| m(A)\| \le \lambda(A)$,  $m$ is a dominated representing  measure (\cite{D}, \cite[p. 148]{BrL}),  and thus $m$  is strongly bounded. 
As in the proof of \cite[Theorem 9]{BCCam},  for every $A\in \Sigma$, $m(A):X\to c_0$ is  compact, thus $p$-convergent.
By hypothesis, $T$ is $p$-convergent.

By Lusin's theorem, there is a compact $K_0\subset K $ so that $\lambda(K\setminus K_0)< 1/4 $ and  $f_n=r_n x_n|_{K_0}$ is continuous for each $n$. Let $H=[f_n]$ be the closed linear subspace spanned by $(f_n)$ in $C(K_0,X)$ and let $S:H\to C(K,X)$  be the isometric extension operator given by   \cite[Theorem 1]{BCCam}. If $h_n= S(f_n)$ for each $n\in \mathbb{N}$, then $(h_n)$ is  
 in the unit ball of $C(K,X)$. By Lemma \ref{le8} (ii), $(f_n)$ is  weakly $p$-summable  in $C(K_0, X)$, and thus $(h_n)$ is weakly $p$-summable   in $C(K, X)$.
%\ref{le13}
%we can apply Lemma \ref{rightnull}; since  $f_n=r_n x_n|_{K_0}$ is continuous, $(r_n|_{K_0})\subset C(C_0)$ and $(x_n)$ is weakly $p$-summable

Let $\hat{T}$ be the extension of $T$ to $B(\Sigma,X)$.   For each $n\ge 0$ we have 
$$\| \hat{T} (r_n x_n) \|=\sup_m \left|\int_K  \langle x_m^*,  r_n(t) x_n \rangle \, r_m(t) \, d\lambda \right|=|\langle x_n^*, x_n\rangle|=1.$$
Therefore
\begin{align*}
\|  T(h_n) \| & \ge \| \hat{T}(r_nx_n)  \| - \| \hat{T}(h_n-r_nx_n)  \| \\
& \ge 1- \| (\int_{K\setminus K_0} \langle x_m^*, h_n(t)-r_n(t)\, x_n \rangle r_m(t)\,d\lambda )_m\| \\
& \ge  1- 2 \lambda(K\setminus K_0) >  1/2.
\end{align*}
Hence $\|T(h_n)\|\not \to 0$, and thus $T$ is not $p$-convergent. This contradiction concludes the proof. 
\end{proof}

A Banach space $X$ has the \emph{Gelfand-Phillips (GP) property} (or $X$ is a \emph{Gelfand-Phillips space}) if every limited subset of $X$ is relatively compact.

Schur spaces and separable spaces have the Gelfand-Phillips  property \cite{BD}.

Let $1\le p<\infty$. A Banach space $X$ has the $p$-\emph{Gelfand-Phillips ($p$-$GP$)  property} (or is a $p$-\emph{Gelfand-Phillips space}) if every  limited weakly $p$-summable  sequence in $X$ is norm null \cite{FZ2}. 

If $X$ has the $GP$ property, then $X$ has the $p$-$GP$ property for any $1\le p<\infty$.

%related to Prop 3.4, Castillo and Sanchez

\begin{proposition} \label{prop15}
Let $K$ be a compact Hausdorff space and $Y$ a Banach space. 

(i) Let $1\le p\le \infty$. If $X\in C_p$, then every strongly bounded operator $T: C(K,X)\to Y$ is $p$-convergent. 

(ii) Let $1\le p< \infty$. If $X$ has the $GP$ property (resp. the $p$-$GP$ property),  then every strongly bounded operator $T: C(K,X)\to Y$ is limited completely continuous (resp. limited $p$-convergent).
\end{proposition}

\begin{proof}  
%The proof is contained in the proof of \cite[Proposition 3.4]{CS}. We include a proof for the sake of completeness.
(i)   Let $(f_n)$ be a weakly $p$-summable sequence in $C(K, X)$. Then for each $t\in K$, $(f_n(t))$ is weakly $p$-summable in $X$, thus norm null. Since $T$ is strongly bounded and $\|f_n(t)\|\to 0$ for each $t\in K$,  $ \|T(f_n)\|\to 0 $ by  \cite[Theorem 2.1]{BrLe}. Thus $T$ is $p$-convergent. 

(ii) Assume that $X$ has the $GP$ property and $T: C(K,X)\to Y$ is strongly bounded. Let $(f_n)$ be a limited weakly null sequence in $C(K, X)$. Then for each $t\in K$, $(f_n(t))$ is limited weakly null in $X$ (by Lemma \ref{lemmalcc}), and thus norm null. 
 %(by the proof of Theorem \ref{Tlcc})
Since $T$ is strongly bounded and $\|f_n(t)\|\to 0$ for each $t\in K$,  $ \|T(f_n)\|\to 0 $ by  \cite[Theorem 2.1]{BrLe}. Thus $T$ is limited completely continuous. 
\end{proof} 

%$1< p<\infty$
\begin{corollary} (\cite[Proposition 2.4]{CS})
Let $1\le p\le \infty$. The following are equivalent about a Banach space $X$:

(a) $X\in C_p$.

(b) For any compact Hausdorff space $K$ and any Banach space $Y$, an operator $m\leftrightarrow T: C(K,X)\to Y$ is $p$-convergent if 
and only if  (i) $T$ is strongly bounded and (ii) for each $A\in \Sigma$, $m(A):X\to Y$ is $p$-convergent.
\end{corollary}

\begin{proof} 
$(a) \Rightarrow (b)$ If $m\leftrightarrow T: C(K,X)\to Y$ is $p$-convergent, then (i) and (ii) follow 
by  \cite[Proposition 2.1]{CS}.  The converse follows from Proposition \ref{prop15}. 
%Conversely, suppose  that (i) and (ii) hold. 

$(b) \Rightarrow (a)$ Suppose that (b) holds and  $X\not \in C_p$. 
Let $K=\Delta=\{-1, 1\}^\mathbb{N}$ be the Cantor set with $\lambda $ be the Haar measure on $\Delta$ (as in the proof of \cite[Theorem 2.1]{PSaab}). Since $K$ is uncountable and metrizable (\cite{L}, p. 35, 36),  $K$ is not dispersed (\cite{Se}, 8.5.5).  
By the proof of Theorem \ref{Th14}, $(iii)  \Rightarrow (i)$, there exists an operator $T:C(K,X) \to c_0$ such that $T$ is strongly bounded, 
 $m(A):X\to c_0$ is $p$-convergent for each $A\in \Sigma$, and $T$ is not $p$-convergent. This contradiction with (b) concludes the proof. 
\end{proof}

If  $m\leftrightarrow T:C(K, X)\to Y$ is completely continuous (resp. $p$-convergent), then $T$ is strongly bounded and for each $A\in \Sigma$, $m(A):X\to Y$ is completely continuous (resp. $p$-convergent) \cite{ID} (resp. \cite[Proposition 2.1]{CS}). 

\begin{corollary} 
Let $1< p<\infty$ and let $K$ be a compact Hausdorff space. The following statements are equivalent:

(a) $K$ is dispersed.

%$DP^*P$ (resp. the $DP^*P_p$) (i) (ii) such  that  $X$ does not have the Schur property

(b) There exists a Banach space $X$ that does not have the Schur property (resp. $X\not \in C_p$) such that  a strongly bounded operator $m\leftrightarrow T:C(K, X)\to c_0$ is weak$^*$ Dunford-Pettis (resp. weak$^*$ $p$-convergent) if and only if for each $A\in \Sigma$, $m(A):X\to c_0$ is weak$^*$ Dunford-Pettis (resp. weak$^*$ $p$-convergent).
\end{corollary}

\begin{proof} 
$(a)\Rightarrow (b)$ by Corollary \ref{Cem1}. 
%{mwDP}
$(b)\Rightarrow (a)$ 
%Suppose (b) holds.
 Suppose $T$ is strongly bounded and for each $A\in \Sigma$, $m(A):X\to c_0$ is completely continuous (resp. $p$-convergent). Since for each $A\in \Sigma$, $m(A):X\to c_0$ is weak$^*$ Dunford-Pettis (resp. weak$^*$ $p$-convergent), 
 $T:C(K, X)\to c_0$ is weak$^*$ Dunford-Pettis (resp. weak$^*$ $p$-convergent) by hypothesis. 
Hence $T$ is completely continuous (resp. $p$-convergent) \cite[Corollary 10]{GRen} (resp. \cite[Corollary 23]{GRen}). Thus an operator $m\leftrightarrow T:C(K, X)\to c_0$ is completely continuous (resp. $p$-convergent) if and only if  $T$ is strongly bounded and for each $A\in \Sigma$, $m(A):X\to c_0$ is completely continuous.   Since $X$ does not have the Schur property (resp. $X\not \in C_p$), 
it follows that $K$ is dispersed by \cite[Theorem 11]{BCCam} (resp. Theorem \ref{Th14}). 
\end{proof}

%In \cite{} it was

%\newpage

\end{document}